\documentclass[10pt,centertags,oneside]{amsart}
\usepackage{amssymb}

\usepackage{times}
\usepackage[mathcal]{euscript}

\copyrightinfo{2002}{V.~Kostrykin, K.A.~Makarov, A.K.~Motovilov}

\renewcommand{\epsilon}{\varepsilon}

\newcommand{\filnam}[1]{\mbox{\texttt{\ignorespaces#1\unskip}}}

\numberwithin{equation}{section}

\newfont{\timitfont}{pplbi7t scaled 1000}
\renewcommand{\v}{\mbox{\timitfont{v}}}

\newcommand{\N}{\mathbb{N}}
\newcommand{\R}{\mathbb{R}}
\newcommand{\EE}{\mathsf{E}}
\newcommand{\bA}{{\mathbf A}}
\newcommand{\bB}{{\mathbf B}}
\newcommand{\bV}{{\mathbf V}}
\newcommand{\cB}{{\mathcal B}}
\newcommand{\cG}{{\mathcal G}}
\newcommand{\cM}{{\mathcal M}}
\newcommand{\cN}{{\mathcal N}}
\newcommand{\cR}{{\mathcal R}}
\newcommand{\cS}{{\mathcal S}}
\newcommand{\cT}{{\mathcal T}}
\newcommand{\cX}{{\mathcal X}}
\newcommand{\cW}{{\mathcal W}}
\renewcommand{\Im}{{\ensuremath{\mathrm{Im}}}}
\newcommand{\fG}{\mathfrak{G}}
\newcommand{\fH}{\mathfrak{H}}
\newcommand{\fK}{\mathfrak{K}}
\newcommand{\fL}{\mathfrak{L}}
\newcommand{\fM}{\mathfrak{M}}
\newcommand{\fN}{\mathfrak{N}}
\newcommand{\fT}{\mathfrak{T}}
\DeclareMathOperator{\Ran}{\mathrm{Ran}}
\DeclareMathOperator{\Ker}{\mathrm{Ker}}
\DeclareMathOperator{\Dom}{\mathrm{Dom}}
\newcommand{\spec}{{\ensuremath{\mathrm{spec}}}}

\newtheorem{theorem}{Theorem}[section]
\newtheorem{lemma}[theorem]{Lemma}
\newtheorem{proposition}[theorem]{Proposition}
\newtheorem{corollary}[theorem]{Corollary}
\newtheorem{remark}[theorem]{Remark}
\newtheorem{hypothesis}[theorem]{Hypothesis}
\newtheorem{definition}[theorem]{Definition}
\newtheorem{example}[theorem]{Example}{\bf}{\rm}

\title[Operator Riccati Equation. A Geometric Approach]
{Existence and Uniqueness of Solutions to the Operator Riccati Equation. A Geometric Approach}

\author[V.~Kostrykin]{Vadim Kostrykin}
\address{Vadim Kostrykin,
Fraunhofer-Institut f\"{u}r Lasertechnik, Steinbachstra{\ss}e
15, D-52074\\ Aachen, Germany}
\email{kostrykin@t-online.de, kostrykin@ilt.fraunhofer.de}

\author[K.~A.~Makarov]{Konstantin A.~Makarov}
\address{Konstantin A.~Makarov,
Department of Mathematics, University of Missouri, Co\-lum\-bia, MO 65211, USA}
\email{makarov@math.missouri.edu}

\author[A.~K.~Motovilov]{Alexander K.~Motovilov}
\address{Alexander K.~Motovilov,
Joint Institute for Nuclear Research, 141980 Dubna, Moscow Region, Russia}
\curraddr{Department of Mathematics, University of Missouri, Columbia, MO 65211, USA}
\email{motovilv@thsun1.jinr.ru}

\keywords{Perturbation theory, spectral subspaces,
graph subspaces, operator Riccati equation}

\subjclass{Primary 47A15, 47A55, 47A62; Secondary 47A53}

\date{July 15 2002}

\begin{document}

\begin{abstract}
We introduce a new concept of unbounded solutions to the operator Riccati equation
$A_1X-XA_0-XVX+V^\ast=0$ and give a complete description of its solutions associated with the
spectral graph subspaces of the block operator matrix $\bB=\begin{pmatrix} A_0 & V \\ V^\ast & A_1
\end{pmatrix}$. We also provide a new characterization of the set of all contractive solutions
under the assumption that the Riccati equation has a contractive solution associated with a
spectral subspace of the operator $\bB$. In this case we establish a criterion for the uniqueness
of contractive solutions.
\end{abstract}

\maketitle

\section{Introduction}\label{sec:1}

In the present article we address the problem of a perturbation of invariant subspaces of
self-adjoint operators on a separable Hilbert space $\fH$ and related questions of the existence
and uniqueness of solutions to the operator Riccati equation.

Given a self-adjoint operator $\bA$ and its closed invariant subspace
$\fH_0\subset\fH$ we set $A_i=\bA|_{\fH_i}$, $i=0,1$ with $\fH_1=\fH\ominus\fH_0$.
Assuming that the perturbation $\bV$ is off-diagonal with respect to the orthogonal
decomposition $\fH=\fH_0 \oplus \fH_1$ consider the self-adjoint operator
\begin{equation*}
\bB = \bA + \bV =\begin{pmatrix} A_0 & V \\ V^\ast & A_1 \end{pmatrix}
\end{equation*}
where $V$ is a linear operator from $\fH_1$ to $\fH_0$. It is well known (see, e.g.,
\cite{Albeverio}, \cite{Apostol:Foias:Salinas}, \cite{Daughtry}) that the Riccati equation
\begin{equation}\label{intro:Ricca}
A_1 X - X A_0 -X V X + V^\ast =0
\end{equation}
has a bounded solution $X:\fH_0 \rightarrow\fH_1$ iff its graph
\begin{equation}\label{intro:Graph}
\cG(\fH_0,X) := \left\{x_0\oplus X x_0|\ x_0\in\fH_0 \right\}
\end{equation}
is an invariant subspace for the operator $\bB$. It might happen, however, that the operator $\bB$
has invariant subspaces that are the graphs of closed densely defined unbounded operators $X:\fH_0
\rightarrow\fH_1$, and the problem of more general solutions to the Riccati equation naturally
arises.

In the present article we introduce the new concept of \emph{unbounded} (closed
densely defined) operator solutions to the Riccati equation and we obtain a geometric
criterion for their existence (Corollary \ref{Riccati:main}) resulting in the
complete description of the bijective correspondence between solutions of the Riccati
equation and the $\bB$-invariant graph subspaces.

Among all solutions to the Riccati equation, those corresponding
to the spectral subspaces of the operator $\bB$, i.e., the
solutions $X$ such that $\cG(\fH_0, X)=\Ran\EE_{\bB}(\Delta)$, the
range of the spectral projection corresponding to some Borel set
$\Delta\subset\R$\,, are of particular interest. Using the
Douglas-Pearcy theorem \cite{Douglas:Pearcy} we prove that a
solution to the Riccati equation is associated with a spectral
subspace iff it is an isolated point in the set of all its
solutions (Theorem \ref{uni}).

Revisiting the case of bounded solutions, we give a complete description of the set of all
contractive solutions ($\|X\|\leq 1$) to the Riccati equation, provided that the Riccati equation
has a contractive solution which is associated with a spectral subspace (Lemma \ref{Neu:Theorem}
and Theorem \ref{last:but:not:least}). This result substantially generalizes the recent uniqueness
theorem due to Adamyan, Langer, and Tretter \cite{Adamyan:Langer:Tretter:2000a}.

In the forthcoming publications
\cite{alpha,Kostrykin:Makarov:Motovilov:4,Kostrykin:Makarov:Motovilov:3} we prove a
number of new existence and uniqueness results for solutions of the Riccati equation
assuming some conditions on the spectra of the operators $A_0$ and $A_1$. Also we
obtain sharp estimates for the norm of these solutions. These estimates are related
to the study of the subspace perturbation problem
\cite{Kostrykin:Makarov:Motovilov:1}.

To avoid getting into technical issues that may obscure the basic ideas of the work,
we assume in this paper that the operator $\bA$ and the perturbation $\bV$ are
bounded. In some cases this hypothesis can easily be relaxed to handle the case of
unbounded $\bA$'s and even unbounded perturbations $\bV$ as well. The extension to
unbounded operators will be presented elsewhere.

The article is organized as follows. In Section \ref{sec:2} we collect some known facts about two
closed subspaces of a separable Hilbert space $\fH$. In Section \ref{sec:3} we give a particularly
simple proof of the Halmos theorem \cite{Halmos:69} providing a criterion for a closed subspace of
the Hilbert space $\fH$ to be the graph $\cG(\fH_0,X)$ of a closed densely defined operator $X$
from a closed subspace $\fH_0\subset \fH$ to its orthogonal complement $\fH_1=\fH_0^\perp$. In
Section \ref{sec:4} we formulate and prove a general criterion for the solvability of the operator
Riccati equation in the class of closed densely defined (not necessarily bounded) operators. The
structure of the set of all solutions to the Riccati equation is analyzed in Section \ref{sec:5}
from the topological point of view. Finally, Section \ref{sec:6} is devoted to the thorough
analysis of the set of all contractive solutions going beyond the one undertaken recently in
\cite{Adamyan:Langer:Tretter:2000a}. In particular, we establish a general criterion for a
contractive solution which is associated with a spectral subspace of the operator matrix $\bB$ to
be unique with no additional assumptions on the spectra of the operators $A_0$ and $A_1$.

A few words about the notations used throughout the paper. Given a linear operator $A$ on a
Hilbert space $\fK$, by $\spec(A)$ we denote the spectrum of $A$. If not explicitly stated
otherwise, $\fN^\perp$ denotes the orthogonal complement in $\fK$ of a subspace $\fN\subset\fK$,
i.e., $\fN^\perp=\fK\ominus\fN$. The identity operator on $\fK$ is denote by $I_\fK$. The notation
$\cB(\fK,\fL)$ is used for for the Banach algebra of bounded operators from the Hilbert space
$\fK$ to the Hilbert space $\fL$. Finally, we write $\cB(\fK)=\cB(\fK,\fK)$.

\subsection*{Acknowledgments.}
V.~Kostrykin is grateful to V.~Enss, A.~Knauf, H.~Leschke, and R. Schrader for useful discussions.
A.~K.~Motovilov acknowledges the great hospitality and financial support by the Department of
Mathematics, University of Missouri--Columbia, USA. He was also supported in part by the Russian
Foundation for Basic Research within Project RFBR 01-01-00958.

\section{Geometry of Two Subspaces of the Hilbert Space}\label{sec:2}
\setcounter{equation}{0}

In this section we collect some facts about pairs of closed subspaces of a separable Hilbert
space. Although most of them are well known, they are scattered in the literature and frequently
formulated in a different form which does not fit the context of the present paper. Without any
attempt to give a complete overview of the whole work done in this direction we mention the
pioneering work of Friedrichs \cite{Friedrichs}, M.~Krein, Krasnoselsky, and Milman
\cite{Krein:Krasnoselsky}, \cite{Krein:Krasnoselsky:Milman}, Dixmier \cite{Dixmier:48},
\cite{Dixmier:49}, Davis \cite{Davis:58}, and Halmos \cite{Halmos:69}. Some of the results
described in this section admit an extension to the case of Banach spaces. We refer the interested
reader to the papers \cite{Krein:Krasnoselsky:Milman} and \cite{Kalton}.
\begin{definition}\label{projections}
Let $(P,Q)$ be an ordered pair of orthogonal projections in $\fH$. We use the
standard notation as introduced by Halmos \cite{Halmos:69} (see also
\cite{Spitkovsky})
\begin{align*}
\fM_{pq}& :=\left\{f\in\fH\big| Pf=pf,\ Qf=qf \right\},\quad p,q=0,1
\\
\fM'_0& :=\Ran P\ \ominus\
(\fM_{10}\ \oplus\ \fM_{11})
\\
\fM'_1& :=\Ran P^\perp\ \ominus\
(\fM_{00}\ \oplus\ \fM_{01})
\\
\fM'& :=\fM'_0\oplus\fM'_1
\\
P'& :=P|_{\fM'}
\\
Q'& :=Q|_{\fM'}.
\end{align*}
\end{definition}

To avoid possible confusion we will often write $\fM_{pq}(P,Q)$ instead of the
shorthand notation $\fM_{pq}$ to emphasize that the canonical decomposition of the
Hilbert space $\fH$ is considered with respect to the ordered pair $(P,Q)$.

{}Following Halmos \cite{Halmos:69} we call the pair $(P',Q')$ the generic part of the pair
$(P,Q)$. Roughly speaking, $(P',Q')$ is the non-commuting part of $(P,Q)$.  Indeed, if $P$ and $Q$
commute, then $P'=Q'=0$.

\begin{theorem}\label{2.2}
Let $(P,Q)$ be an ordered pair of orthogonal projections in the
Hilbert space $\fH$. Then the space $\fH$ admits the
(canonical) orthogonal decomposition
\begin{equation}\label{CanonDecomp}
\fH=\fM_{00}\ \oplus\ \fM_{01}\ \oplus\ \fM_{10}\ \oplus\
\fM_{11}\ \oplus\ \fM'.
\end{equation}
With respect to this decomposition the projections $P$ and $Q$ read
{\setlength\arraycolsep{2pt}
\begin{equation*}
\begin{array}{ccccccccccc}
P &=& 0 & \oplus & 0            & \oplus & I_{\fM_{10}} & \oplus &
I_{\fM_{11}} & \oplus & P',\\
Q &=& 0 & \oplus & I_{\fM_{01}} & \oplus & 0            & \oplus &
I_{\fM_{11}} & \oplus & Q'.
\end{array}
\end{equation*}}
With respect to the decomposition $\fM'=\fM'_0 \oplus \fM'_1$ the projections $P'$
and $Q'$ read
\begin{equation}\label{generic:part}
P'\quad =\quad \begin{pmatrix} I_{\fM'_0} & 0
\\ 0 & 0
\end{pmatrix}, \qquad
Q'\quad =\quad \cW^\ast \begin{pmatrix} \cos^2\Theta & \sin\Theta \cos\Theta
\\ \sin\Theta \cos\Theta & \sin^2\Theta
\end{pmatrix} \cW ,
\end{equation}
where $\Theta$ is a (unique) positive semidefinite angle operator in $\fM'_0$
such that
\begin{equation*}
\sin^2 \Theta= P'(I_{\fM'}-Q') P'|_{\fM'_0};
\end{equation*}
$\spec(\Theta)\subset[0,\pi/2]$ but $0$ and $\pi/2$ are not eigenvalues of $\Theta$.
The unitary operator matrix $\cW:\ \fM'_0\oplus\fM'_1\rightarrow \fM'_0\oplus\fM'_0$
reads
\begin{equation*}
\cW=\begin{pmatrix} I_{\fM'_0} & 0 \\ 0 & W
\end{pmatrix}
\end{equation*}
where $W\in\cB(\fM_1',\fM_0')$ is the unitary operator from
the polar decomposition
\begin{equation*}
{P'}^\perp Q'P'\vert_{\fM_0'}=W^\ast \left(({P'}^\perp Q'P'\vert_{\fM_0'})^\ast {P'}^\perp
Q'P'\vert_{\fM_0'}\right)^{1/2}.
\end{equation*}

In particular, the difference $Q'-P'$ of the generic parts of
the projections $P$ and $Q$ can be represented with respect to
the decomposition $\fM'=\fM'_0\oplus\fM'_1$ in the form
\begin{equation}\label{angle:operator}
\begin{split}
Q'-P' &=\cW^\ast \begin{pmatrix} \sin\Theta & 0
\\ 0 & \sin\Theta
\end{pmatrix} \begin{pmatrix} -\sin\Theta & \cos\Theta
\\ \cos\Theta & \sin\Theta
\end{pmatrix} \cW\\
&=\cW^\ast \begin{pmatrix} -\sin\Theta & \cos\Theta
\\ \cos\Theta & \sin\Theta
\end{pmatrix} \begin{pmatrix} \sin\Theta & 0
\\ 0 & \sin\Theta
\end{pmatrix} \cW.
\end{split}
\end{equation}
and hence
\begin{equation}\label{P-Q}
\|Q'-P'\|=\|\sin
\Theta(Q', P')\|.
\end{equation}
\end{theorem}
In a slightly different form Theorem \ref{2.2} was proven by Davis \cite{Davis:58} and Halmos in
\cite{Halmos:69}. An alternative, simple and direct proof of this theorem was given by Amrein and
Sinha \cite{Amrein:Sinha}.

Theorem \ref{2.2} has been proved to be of great importance in a number of problems related to
pairs of orthogonal projections. In particular, it was successfully used for the study of the
operator algebras generated by a pair of orthogonal projections (see \cite{Spitkovsky},
\cite{Takesaki} and references therein).

In the next section we will study the graph subspaces associated with an orthogonal
decomposition of the Hilbert space and we will revisit Theorem \ref{2.2} which allows
to perform the subsequent analysis in a particularly simple manner.

\section{Graph Subspaces}\label{sec:3}
\setcounter{equation}{0}

\begin{definition}\label{DefGraf}
Let $\fH_0$ be a closed subspace of a Hilbert space $\fH$ and $X$ a closed densely defined
(possibly unbounded) operator from $\fH_0$ to $\fH_1=\fH_0^\perp$ with domain $\Dom (X)$. The
closed linear subspace
\begin{equation*}
\cG(\fH_0, X)=\{x\in\fH|\ x= x_0\oplus X x_0,\ x_0\in\Dom (X)\subset\fH_0\}
\end{equation*}
is called the graph subspace of ${\fH}$ associated with the pair
$(\fH_0,X)$ or, in short, the graph of $X$.
\end{definition}

One easily checks that
\begin{equation}
\label{graph} \cG(\fH_0,X)^\perp=\cG(\fH_0^\perp,-X^\ast).
\end{equation}

We start with presenting a fairly simple and partly known result
(see~\cite{Halmos:69}) that characterizes the graph subspaces in
terms of the canonical decomposition \eqref{CanonDecomp}.
\begin{theorem}\label{3.2}
Let $P$ and $Q$ be orthogonal projections in a Hilbert space
$\fH$. The subspace $\Ran Q$ is a graph subspace $\cG(\Ran P,X)$
associated with some closed densely defined (possibly unbounded)
operator $X: \Ran P\rightarrow \Ran P^\perp$ with
$\Dom(X)\subset\Ran P$ iff the subspaces $\fM_{01}(P,Q)$ and
$\fM_{10}(P,Q)$ in the canonical decomposition of the Hilbert
space $\fH$ \eqref{CanonDecomp} are trivial, i.e.,
\begin{equation}\label{triv}
\fM_{01}(P,Q) = \fM_{10}(P,Q) = \{0\}.
\end{equation}
For given orthogonal projection $P$ the correspondence between
the closed subspaces $\Ran Q$ satisfying \eqref{triv} and closed
densely defined operators $X:\Ran\to\Ran P^\perp$ is one-to-one.
\end{theorem}
\begin{proof}
\textit{``If"} Part. Assume \eqref{triv}. Let $P'$ and $Q'$ be
generic parts of the projections $P$ and $Q$, respectively.
{}From \eqref{generic:part} it follows that $\Ran Q'$ given by
\begin{equation*}
\Ran Q' = \{ x\in\fM'\, |\, \cos \Theta\ x_0 + W^\ast \sin \Theta\ x_0,\ x_0\in\fM'_0\}
\end{equation*}
is a graph subspace of the generic subspace $\fM'=\Ran P'\oplus
(\Ran P')^\perp$.
More explicitly,
\begin{equation*}
\Ran Q'=\cG(\Ran P', W^\ast \tan \Theta)
\end{equation*}
with
\begin{equation*}
\Dom (\tan \Theta)=\{x_0\in \Ran P' \, |\, x_0 = P' x \text{ for some } x\in \Ran Q'\}.
\end{equation*}
Introducing the operator $X$
from $\Ran P$ to $\Ran P^\perp$ with
\begin{equation}
\label{DomX}
\Dom(X)=
\Dom(\tan \Theta) \oplus (\Ran P\ominus \fM_0')
\end{equation}
by
\begin{equation}\label{defX}
X x = \begin{cases} W^\ast \tan \Theta\ x, & x\in \fM'_0, \\ 0, & x\in \Ran P\ominus \fM'_0
\end{cases}
\end{equation}
yields $\Ran Q = \cG(\Ran P, X)$ since \eqref{triv} holds.

\textit{``Only If"} Part. Assume that $\Ran Q$ is a graph
subspace associated with a closed densely defined operator $X$,
i.e., $\Ran Q = \cG(\Ran P, X)$. To prove \eqref{triv} it
suffices to establish that the points $\pm 1$ are not
eigenvalues of $Q-P$, i.e., $\Ker (Q-P\pm I_{\fH})=0$.

Suppose to the contrary that, say, $+1$ is an eigenvalue of $Q-P$, that
is,
\begin{equation}\label{P}
(Q-P)f=f, \quad 0\neq f\in\Ran P,
\end{equation}
and, hence, by \eqref{graph} $f$ admits the decomposition
\begin{equation}\label{Q}
f=x+Xx -X^\ast y+y
\end{equation}
for some $x\in \Dom (X) \subset\fN$ and $y\in \Dom (X^\ast)\subset\Ran P^\perp$. By inspection
\begin{equation}\label{PQ}
(Q-P)(x+Xx -X^\ast y+y)=Xx+X^\ast y.
\end{equation}
Therefore, combining \eqref{P}, \eqref{Q},
and \eqref{PQ} yields
\begin{equation*}
X x+X^\ast y=x+Xx -X^\ast y+y
\end{equation*}
and
\begin{equation*}
2X^\ast y=x+y,
\end{equation*}
which is only possible if $x=y=0$ and, thus, $f=0$.
Hence, the point $+1$ is not an eigenvalue for $Q-P$.

One proves that $-1$ is not an eigenvalue of $Q-P$ in a similar
way.

The last statement of the theorem follows from the fact that if
two closed graph subspaces $\cG(\Ran P,X_1)$ and $\cG(\Ran P,X_2)$
coincide iff $X_1=X_2$ (see, e.g., \cite{Kato}).
\end{proof}
\begin{remark}
Under the hypothesis of Theorem \ref{3.2}
\begin{align}
\fM_{11}(P,Q) &= \Ker X,\label{ker}\\
\fM_{00}(P,Q)& = \Ker X^\ast = (\Ran X)^\perp.\label{koker}
\end{align}
The first representation holds by the definition \eqref{defX} of
the operator $X$ and \eqref{koker} follows from \eqref{ker} by
duality argument \eqref{graph}.
\end{remark}

The following reformulation of Theorem \ref{3.2} distinguishes
the cases of the graph subspaces associated with bounded and
unbounded operators $X$, respectively.

\begin{corollary}\label{frpr}
Assume Hypothesis \ref{projections}. Then:

\item{(i)} The inequality $\|P-Q\|<1$ holds true iff $\Ran Q$ is
a graph subspace associated with the subspace $\Ran P$ and some
bounded operator $X\in\cB(\Ran P,\Ran P^\perp)$, that is, $\Ran
Q=\cG(\Ran P, X)$. In this case

\begin{equation}
\label{glgl}
\|X\| = \frac{\|P-Q\|}{\sqrt{1-\|P-Q\|^2}}
\end{equation}
and
\begin{equation}
\label{estimate}
\|P-Q\|=\frac{\|X\|}{\sqrt{1+\|X\|^2}}.
\end{equation}

\item{(ii)} $\fM_{10}(P,Q) = \fM_{01}(P,Q) = \{0\}$ and
$\|P-Q\|=1$ iff $\Ran Q$ is a graph subspace associated with the
subspace $\Ran P$ and an unbounded operator $X$ from $\Ran P$ to
$\Ran P^\perp$, i.e., $\Ran Q=\cG(\Ran P, X)$.
\end{corollary}
\begin{proof} (i). By Theorem \ref{3.2} $\Ran Q$ is a graph
subspace with respect to the projection $P$ if and only if
$\fM_{10}(P,Q) = \fM_{01}(P,Q) = \{0\}$, and hence
\begin{equation*}
\|P-Q\|=\|P'-Q'\|,
\end{equation*}
where $(P',Q')$ is the generic part of the pair $(P,Q)$. By Theorem
\ref{2.2}
\begin{equation*}
\|P'-Q'\|=\|\sin \Theta(P',Q')\|,
\end{equation*}
where $\Theta(P',Q')$ is the operator angle between the
subspaces $\Ran P'$ and $\Ran Q'$ and, moreover,
\begin{equation*}
\Ran Q' = \cG(\Ran P', W^\ast \tan \Theta(P',Q')), \qquad \Ran Q = \cG(\Ran P, X),
\end{equation*}
where $X$ is the extension of $W^\ast \tan \Theta(P',Q')$ given by \eqref{DomX}, \eqref{defX}.

Clearly, $X$ is bounded iff the operator $\tan \Theta(P', Q')$
is bounded. The equality \eqref{defX} implies
\begin{equation}\label{tri}
\|X\|=\|\tan \Theta(P',Q')\|,
\end{equation}
and then
\eqref{estimate} is a consequence of the trigonometric identity
\begin{equation*}
\sin \theta=\frac{\tan\theta}{\sqrt{1 + \tan^2\theta}},\quad \theta \in
[0, {\pi}/{2})
\end{equation*}
combining \eqref{raz}, \eqref{dva}, and \eqref{tri} which proves (i).

(ii). The operator $X$ is unbounded iff $\pi/2\in \spec
(\Theta(P',Q'))$. In this case $\|P'-Q'\|=\|P-Q\|=1$ by
\eqref{dva} and \eqref{tri} which proves (ii).
\end{proof}

\begin{remark}\label{folklore}
Part (i) with inequality sign instead of the equality \eqref{estimate} is well known.  A proof can
be found, e.g., in \cite[Theorem 1]{Daughtry} or \cite[Lemma 2.3]{Apostol:Foias:Salinas}.
\end{remark}
\begin{remark}
The orthogonal projection $Q$ onto the graph subspace $\cG(\Ran
P,X)$ corresponding to a closed densely defined operator
$X:\,\Ran P\rightarrow\Ran P^\perp$ can be written as $2\times
2$ operator matrix with respect to the orthogonal decomposition
$\fH=\Ran P\oplus\Ran P^\perp$
\begin{equation}\label{qu}
Q = \begin{pmatrix}
  (I_{\fH_0} + X^\ast X)^{-1} & \overline{(I_{\fH_0} + X^\ast X)^{-1} X^\ast} \\[1ex]
  X (I_{\fH_0} + X^\ast X)^{-1} & \overline{X (I_{\fH_0} + X^\ast X)^{-1} X^\ast}
\end{pmatrix},
\end{equation}
where $\fH_0=\Ran P$ and the bar denotes the closure. The operator entries of \eqref{qu} are
bounded operators since $\Dom(X^\ast X)\subset \Dom(X)$ is a core for $X$ (see, e.g.,
\cite{Kato}).
\end{remark}

\section{Riccati Equation}\label{sec:4}
\setcounter{equation}{0}

The main purpose of this section is to introduce a concept of
closed densely defined (possibly unbounded) operator solutions to
the Riccati equation and to provide a geometric criterion of their
existence.

Throughout this section we adopt the following hypothesis.
\begin{hypothesis}\label{diag}
Assume that the separable Hilbert space $\fH$ is decomposed into
the orthogonal sum of two subspaces
\begin{equation}\label{decom}
\fH=\fH_0 \oplus \fH_1.
\end{equation}
Assume, in addition, that $\bB$ is a self-adjoint operator
represented with respect to the decomposition \eqref{decom} as a
$2\times2$ operator block matrix
\begin{equation*}
\bB=\begin{pmatrix}
  A_0 & V \\
  V^\ast & A_1
\end{pmatrix},
\end{equation*}
where $A_i\in\cB(\fH_i)$, $i=0,1$, are bounded self-adjoint operators in $\fH_i$ while
$V\in\cB(\fH_1,\fH_0)$ is a bounded operator from $\fH_1$ to $\fH_0$. More explicitly, $\bB
=\bA+\bV$, where $\bA$ is the bounded diagonal self-adjoint operator,
\begin{equation*}
\bA = \begin{pmatrix}
  A_0 & 0 \\
  0 & A_1
\end{pmatrix},
\end{equation*}
and the operator $\bV=\bV^\ast$ is an off-diagonal bounded operator
\begin{equation*}
\bV=\begin{pmatrix}
  0 & V \\
  V^\ast & 0
\end{pmatrix}.
\end{equation*}
\end{hypothesis}

The notion of strong and weak \emph{bounded} solutions to the Riccati equation with
\emph{unbounded} operator coefficients was introduced in \cite{Albeverio} (cf.~\cite{Lasiecka}). In our case where the
operator coefficients are bounded but solutions are allowed to be unbounded we use
the following definition.

\begin{definition}
\label{DefRicc}
Assume Hypothesis \ref{diag}.
A closed densely defined (possibly unbounded) operator $X$ from $\fH_0$
to $\fH_1$ with $\Dom (X)$ is called a weak solution to the
Riccati equation
\begin{equation}\label{Riccati}
A_1 X - X A_0 - XVX + V^\ast =0
\end{equation}
if for any $x\in\Dom (X)$ and any $y\in\Dom (X^\ast)$
\begin{equation}
\label{Ricform} (A_1 y,\ X x) - (X^\ast y,\ A_0 x) - (X^\ast y, VX x) + (Vy, x) = 0.
\end{equation}

A closed densely defined (possibly unbounded) operator
$X$ from $\fH_0$ to $\fH_1$ with $\Dom (X)$ is called a strong solution to the
Riccati equation \eqref{Riccati} if
\begin{equation*}
\left.\Ran(A_0+VX)\right|_{\Dom(X)}\subset\Dom(X)
\end{equation*}
and
\begin{equation*}
A_1 X x - X (A_0 + VX) x + V^\ast x =0\quad\text{for any}\quad x\in\Dom(X).
\end{equation*}
\end{definition}

Obviously, if $X$ is a bounded operator, then the Riccati equation
\eqref{Riccati} can be understood as an operator equality.

The notions of weak and strong solutions to the Riccati equation are in fact
equivalent. The precise statement is as follows.

\begin{lemma}\label{weak:strong}
Assume Hypothesis \ref{diag}. A closed densely defined (possibly unbounded) operator
$X$ from $\fH_0$ to $\fH_1$ with $\Dom (X)$ is a weak solution to the Riccati
equation \eqref{Riccati} iff $(A_0 x + V X x)\in\Dom(X)$, $x\in\Dom(X)$, and
\begin{equation}
\label{Ric:strong} A_1 X x - X (A_0 + VX) x + V^\ast x =0\quad\text{for any}\quad x\in\Dom(X),
\end{equation}
i.e., $X$ is a strong solution to \eqref{Riccati}.
\end{lemma}

\begin{proof}
Assume that $X$ is a weak solution to the Riccati equation \eqref{Riccati}, i.e., \eqref{Ricform}
holds for any $x\in\Dom (X)$ and any $y\in\Dom (X^\ast)$. Then
\begin{equation*}
(y,\ A_1 X x + V^\ast x)=(X^\ast y,\ A_0 x+V X x),
\end{equation*}
which implies, in particular, that $A_0 x+V X x \in \Dom (X^{\ast\ast})$. Since $X$ is closed and
densely defined, one infers $X^{\ast\ast} = X$ and, therefore, $A_1 X x + V^\ast x = X (A_0+V X)x$
for all $x\in\Dom(X)$.

The converse statement is obvious.
\end{proof}

As a consequence of Lemma \eqref{weak:strong} we obtain the following theorem.

\begin{theorem}\label{2x2}
Assume Hypothesis \ref{diag}. A closed densely defined (possibly unbounded) operator
$X$ from $\fH_0$ to $\fH_1$ with $\Dom (X)$ is a weak solution to the Riccati
equation \eqref{Riccati} iff the graph subspace $\cG(\fH_0, X)$ is invariant for the
operator $\bB$.
\end{theorem}
\begin{proof}
First, assume that $\cG(\fH_0, X)$ is invariant for $\bB$. Then
\begin{equation*}
\bB(x\oplus Xx)=(A_0 x + V X x)\oplus(A_1 X x + V^\ast x)\in\cG(\fH_0, X) \text{ for any } x\in
\Dom (X).
\end{equation*}
In particular, $A_0 x + V X x\in \Dom (X)$ and
\begin{equation*}
A_1 X x + V^\ast x = X (A_0 x + V X x) \text{ for all } x\in \Dom (X).
\end{equation*}
Hence,
\begin{equation*}
(y, V^\ast x+A_1 X x)=(y, X(A_0 x+V X x))\text{ for all } x\in \Dom (X) \text{ and for all } y\in
\Dom (X^\ast),
\end{equation*}
which proves that $X$ is a weak solution to the Riccati equation
\eqref{Riccati}.

To prove the converse statement assume that $X$ is a weak solution to the Riccati equation
\eqref{Riccati}, i.e., \eqref{Ricform} holds for any $x\in\Dom (X)$ and any $y\in\Dom (X^\ast)$.
{}From Lemma \ref{weak:strong} it follows that
\begin{equation*}
A_0 x + V X x \in \Dom (X)
\end{equation*}
and
\begin{equation*}
A_1 X x + V^\ast x = X (A_0 x+ V X x),\qquad x\in \Dom (X),
\end{equation*}
which proves that the graph subspace $\cG(\fH_0, X)$ is $\bB$-invariant.
\end{proof}

The next statement is an immediate corollary of Theorems \ref{frpr} and \ref{2x2}.

\begin{corollary}\label{Riccati:main}
Assume Hypothesis \ref{diag}. Let $\fG$ be a
closed $\bB$-invariant subspace of the Hilbert
space $\fH$ and $P$ and $Q$ denote the orthogonal
projections in $\fH$ respectively onto $\fH_0$ and
$\fG$. Then:

(i) The inequality
\begin{equation*}
\|P-Q\|<1
\end{equation*}
holds iff $\fG$ is a graph subspace,
$\fG=\cG(\fH_0,X)$ where $X$ is a bounded solution
to the Riccati equation \eqref{Riccati}. In this
case equalities \eqref{glgl} and \eqref{estimate}
hold true.

(ii) The equality
\begin{equation*}
\|P-Q\|=1
\end{equation*}
holds and
\begin{equation*}
\fM_{01}(P,Q)=\fM_{10}(P,Q)=\{0\},
\end{equation*}
iff $\fG$ is a graph subspace, $\fG=\cG(\fH_0,X)$,
where X is a closed densely defined unbounded weak
solution to the Riccati equation \eqref{Riccati}.
\end{corollary}

We present an example where the Riccati equation has an unbounded solution.

\begin{example}
Let $\fH=\fH_0 \oplus \fH_1$ where
$\fH_0=\fH_1=L^2(0,1)$. Let $\Lambda$ be the
multiplication operator in $L^2(0,1)$,
\begin{equation*}
(\Lambda f)(\lambda)= \lambda f(\lambda), \quad f\in L^2(0,1),
\end{equation*}
and $A_0=-\Lambda$, $A_1=\Lambda$, and $V=\Lambda^2$. In this case
the Riccati equation \eqref{Riccati} being of the form
\begin{equation*}
\Lambda X + X\Lambda - X\Lambda^2 X +\Lambda^2=0
\end{equation*}
has a unbounded self-adjoint solution $X=f(\Lambda)$ where
\begin{equation*}
f(\lambda) = -\frac{1+\sqrt{1+\lambda^2}}{\lambda}.
\end{equation*}
\end{example}

\section{Solutions Associated With Spectral Subspaces}\label{sec:5}

The structure of the set of solutions to the
Riccati equation associated with spectral
subspaces of the operator $\bB$ can be studied
based on the Douglas-Pearcy theorem \cite[Theorem
3]{Douglas:Pearcy} on invariant subspaces of
normal operators.

\begin{theorem}\label{Douglas:Pearcy}
Let $T$ be a bounded self-adjoint operator in a Hilbert space $\fH$ and $Q$ an
orthogonal projection onto a closed $T$-invariant subspace of $\fH$. Then the
following are equivalent:

(i) $\Ran Q$ is a spectral subspace of the
operator $T$, i.e., there is a Borel set
$\Delta\subset\R$ such that $Q=\EE_T(\Delta)$,
where $\EE_T(\Delta)$ denotes the spectral
projection of $T$ corresponding to the set
$\Delta$;

(ii) $\|Q-P\|=1$ for any orthogonal projection $P$
in $\fH$, $P\neq Q$, such that $\Ran P$ is
$T$-invariant;

(iii) $\dim\fM_{10}(P,Q)+\dim\fM_{01}(P,Q)>0$ for
any orthogonal projection $P\neq Q$ in $\fH$ such
that $\Ran P$ is $T$-invariant;

(iv) $Q$ is an isolated point (in the operator
norm topology) of the set of all orthogonal
projections onto all $T$-invariant subspaces.
\end{theorem}

\begin{proof}
The equivalence of (i), (ii), and (iv) is proven
in \cite{Douglas:Pearcy}. The implication
(iii)$\Rightarrow$(ii) is implied by the
decomposition \eqref{CanonDecomp}. Thus, we will
only prove the implication (i)$\Rightarrow$(iii).

Assume that (i) holds. Suppose to the contrary that (iii) does not hold. That is, $Q$
is a spectral projection for $T$ such that $\fM_{10}(P,Q)=\fM_{01}(P,Q)=\{0\}$ for
some orthogonal projection $P\neq Q$ such that $\Ran P$ is $T$-invariant. Since $T$
is self-adjoint, the subspace $\Ran P$ is reducing for $T$ and thus $TP=PT$.
Therefore, (see, e.g., \cite[Theorem 6.3.2]{Birman:Solomyak}) $P$ commutes with all
spectral projections of $T$. In particular, $PQ=QP$. Since
$\fM_{10}(P,Q)=\fM_{01}(P,Q)=\{0\}$, we conclude that $P=Q$, a contradiction. Thus,
(i) implies (iii).
\end{proof}

In the following we will use a concept of the
generalized convergence of closed operators (see
\cite[Section IV.2]{Kato}). The generalized
convergence is a natural extension to the case of
unbounded operators of a notion of uniform
convergence. For convenience of the reader we
recall its definition adapted to the present
context.

\begin{definition}
A sequence $\{X_n\}_{n\in\N}$ of closed densely
defined operators $X_n$ from $\fH_0$ to $\fH_1$
converges in the generalized sense to a closed
operator $X$ from $\fH_0$ to $\fH_1$ if
\begin{equation*}
\lim_{n\rightarrow\infty} \|Q_n - Q\|=0,
\end{equation*}
where $Q_n$ and $Q$ are orthogonal projections
onto the graph subspaces $\cG(\fH_0, X_n)$ and
$\cG(\fH_0, X)$, respectively.
\end{definition}

The following statement characterizes the set of
solutions to the Riccati equation \eqref{Riccati}
associated with spectral subspaces of the operator
$\bB$.

\begin{theorem}\label{uni}
Assume Hypothesis \ref{diag}. Denote by $\cX$ the
set of all (weak) solutions to the Riccati
equation \eqref{Riccati}. Then:

(i) If $X\in \cX$ and the invariant graph subspace
$\cG(\fH_0, X)$ is a spectral subspace of the
operator $\bB$, i.e., $\cG(\fH_0, X)=\Ran
E_{\bB}(\Delta) $ for some Borel set
$\Delta\subset\R$, then $X$ is an isolated point
of the set $\cX$ in the topology of the
generalized convergence of operators.

(ii) If $X\in \cX$ is a bounded operator, then the
invariant graph subspace $\cG(\fH_0, X)$ is a
spectral subspace iff $X$ is an isolated point of
the set $\cX\cap\cB(\fH_0, \fH_1)$ in the
operator norm topology, i.e., there is a
neighborhood of $X$ in $\cB(\fH_0,\fH_1)$ where
the Riccati equation \eqref{Riccati} has no
solutions except $X$.
\end{theorem}
\begin{proof}
(i) Let $\cG(\fH_0, X)$ be a spectral subspace for
$\bB$ and let $Q$ denote the orthogonal projection
in $\fH$ onto $\cG(\fH_0, X)$. Suppose to the
contrary that $X$ is not an isolated solution,
i.e., there is a sequence $\{X_n\}_{n\in\N}$ of
solutions to \eqref{Riccati} such that $X_n\neq
X$, $n\in\N$ and
\begin{equation*}
\lim_{n\to \infty}\|Q_n-Q\|=0,
\end{equation*}
where $Q_n$, $n\in\N$ denote the orthogonal
projections in $\fH$ onto the $\bB$-invariant
graph subspaces $\cG({\fH}_0,X_n)$. By Theorem
\ref{Douglas:Pearcy} this contradicts the
assumption that $\cG(\fH_0, X)$ is a spectral
subspace for $\bB$ which completes the proof of
(i).

(ii) Since by Theorem IV.2.23 in \cite{Kato} the
generalized convergence of bounded operators
implies its uniform convergence, the ``only if"
part follows from (i). Therefore, we only prove
the ``if" part.

Let $X\in\cB(\fH_0,\fH_1)$ be an isolated bounded
solution to \eqref{Riccati}. Suppose that
$\cG(\fH_0, X)$ is not a spectral subspace for
$\bB$. By Theorem \ref{Douglas:Pearcy} (iv) this
implies that there is a sequence of orthogonal
projections $Q_n$, $n\in\N$ such that $\Ran Q_n$
is $\bB$-invariant and
\begin{equation}\label{Konvergenz}
\lim_{n\rightarrow\infty} \|Q_n - Q\| = 0,
\end{equation}
where $Q$ is the orthogonal projection onto
$\cG(\fH_0,X)$. Equation \eqref{Konvergenz} means
that
\begin{equation*}
\|Q_n - Q\| < 1 - \frac{\|X\|}{\sqrt{1+\|X\|^2}}
\end{equation*}
for $n\in \N$ large enough. Therefore,
\begin{equation*}
\|Q_n-P\|\leq \|Q_n-Q\|+\|Q-P\| =
\|Q_n - Q\| + \frac{\|X\|}{\sqrt{1+\|X\|^2}} < 1,
\end{equation*}
$n\in \N$ large enough, where $P$ denotes the
orthogonal projection in $\fH$ onto $\fH_0$. By
Theorem \ref{3.2} for those $n\in \N$, $\Ran Q_n =
\cG(\fH_0, X_n)$ for some $X_n\in\cB(\fH_1,\fH_0)$
where $X_n$ is a solution to \eqref{Riccati} by Theorem \ref{2x2}. Finally, by
Theorem IV.2.23 in \cite{Kato} equality \eqref{Konvergenz} implies
\begin{equation*}
\lim\limits_{n\rightarrow\infty} \| X_n - X\| = 0,
\end{equation*}
which contradicts the assumption that $X$ is an isolated solution.
\end{proof}

\section{Contractive Solutions Associated With Spectral
Subspaces: Uniqueness Criteria}\label{sec:6}

Corollaries \ref{frpr} and \ref{Riccati:main} imply that under Hypothesis \ref{diag} the Riccati
equation \eqref{Riccati} has a contractive solution $X$ iff the subspaces $\fH_0$ and
$\cG(\fH_0,X)$ are invariant for $\bA$ and $\bB=\bA+\bV$, respectively, and the orthogonal
projections $P$ and $Q$ onto these subspaces satisfy $\|P-Q\|\leq \sqrt{2}/2$. It is also known
\cite{Adamyan:Langer:Tretter:2000a} that under the same hypothesis the Riccati equation
\eqref{Riccati} has a contractive solution iff there exists a self-adjoint involution $\mathbf{J}$
in $\fH=\fH_0\oplus \fH_1$ such that
\begin{equation*}
\bB \mathbf{J}=\mathbf{J}\bB
\end{equation*}
and the subspace $\fH_1$ is maximal $\mathbf{J}$-nonnegative, that is, $\fH_1$ is not properly
contained in another $\mathbf{J}$-nonnegative subspace. In principle, these criteria provide
complete although somewhat implicit characterization of the set $\cS$ of all possible contractive
solutions for the Riccati equation \eqref{Riccati}. The main goal of this section is to obtain new
characterization of the set $\cS$ under the assumption that the Riccati equation has at least one
contractive solution associated with a spectral subspace of the operator matrix $\bB$. As a
by-product of this new description we get some uniqueness results generalizing those obtained in
\cite{Adamyan:Langer:Tretter:2000a}.

We start by stating an auxiliary result describing two contractions such that the
orthogonal projections onto their graphs commute.

\begin{lemma}\label{Neu:Theorem}
Assume Hypothesis \ref{diag}. Let
$X,Y\in\cB(\fH_0,\fH_1)$, $\|X\|\leq 1$,
$\|Y\|\leq 1$ be two contractions such that the
orthogonal projections in $\fH$ onto their graphs
$\cG(\fH_0,X)$ and $\cG(\fH_0, Y)$ commute. Then
\begin{equation}\label{relations}
Y|_{\fL} = -X|_{\fL},\qquad Y|_{\fL^\perp} = X|_{\fL^\perp},
\end{equation}
where
\begin{align}\label{i3}
 \fL&=\Ker(I_{\fH_0}+Y^\ast X)=\Ker(I_{\fH_0}+X^\ast Y)
\end{align}
is a subspace of \, $\Ker(I_{\fH_0}-X^\ast X)\cap\Ker(I_{\fH_0}-Y^\ast Y)$ and
\begin{equation*}
\fL^\perp =\fH_0\ominus \fL.
\end{equation*}
Moreover,
\begin{equation}
\label{i4}
\fL=\Ker(X+Y)\ominus\big(\Ker(X)\cap\Ker(Y)\big)
\end{equation}
and
\begin{equation}\label{iadro}
\fL^\perp=\Ker(X-Y).
\end{equation}
\end{lemma}
\begin{proof}
Note that $x\in \Ker(I_{\fH_0} + Y^\ast X)$ means that
\begin{equation*}
\|x\|^2=-(Y^\ast Xx, x)=-(Xx, Yx),
\end{equation*}
which holds if and only if
\begin{equation*}
Yx=-Xx
\end{equation*}
and
\begin{equation*}
x\in\Ker(I_{\fH_0}-X^\ast X) \cap \Ker(I_{\fH_0}-Y^\ast Y)
\end{equation*}
since both $X$ and $Y$ are contractions.
Hence
\begin{equation}\label{Kontraktion:1}
\Ker(I_{\fH_0}+Y^\ast X)=\Ker(X+Y)\cap \Ker (I_{\fH_0}-X^\ast X)\cap \Ker (I_{\fH_0}-Y^\ast Y).
\end{equation}
By symmetry,
\begin{equation*}
\Ker(I_{\fH_0}+Y^\ast X)=\Ker(I_{\fH_0}+X^\ast Y).
\end{equation*}
Therefore, we have proven equalities \eqref{i3},
and the first equality in \eqref{relations}.

Given an arbitrary $x\in \fL^\perp$, one concludes that $(x, y)=0 $ for any $y\in \fL$. By
\eqref{Kontraktion:1} $y\in\Ker(I_{\fH_0}-X^\ast X)$ and, hence,
\begin{equation*}
(x,y)+(Xx,Xy)=0,\quad y\in \fL,
\end{equation*}
which means that
\begin{equation}\label{eq1}
(x\oplus Xx) \perp (y\oplus Xy),\quad y\in \fL.
\end{equation}

By hypothesis the orthogonal projections onto $\cG(\fH_0,X)$ and
$\cG(\fH_0,Y)$ commute. This implies in particular that
\begin{equation*}
\bigl(\cG(\fH_0,X)\ominus\fT\bigr) \perp \bigl(\cG(\fH_0,Y)\ominus \fT\bigr),
\end{equation*}
where $\fT=\cG(\fH_0,X)\cap \cG(\fH_0,Y).$
Introducing the subspace
\begin{equation*}
\cN=P_0 \big (\cG(\fH_0,X)\ominus \fT \big ),
\end{equation*}
where $P_0$ denotes the canonical projection from
$\fH_0\oplus \fH_1$ onto $\fH_0$, one proves by
inspection that
\begin{equation}\label{pio}
\cN\subset \fL.
\end{equation}
In particular, \eqref{eq1} and \eqref{pio} imply that
\begin{equation}\label{zd}
(x\oplus Xx) \perp (y\oplus Xy),\quad\text{ for any }
x\in\fL^\perp\text{ and any } y\in \cN.
\end{equation}
Since
\begin{equation*}
\cG(\fH_0,X)\ominus \fT =\{ y\oplus Xy \, |\, y\in \cN\}
\end{equation*}
and $x\oplus Xx\in \cG(\fH_0,X)$, condition \eqref{zd} means that
\begin{equation*}
x\oplus Xx\in \fT,\quad \text{ for any } x\in\fL^\perp.
\end{equation*}
Therefore, by the definition of the subspace $\fT$,
\begin{equation*}
Xx=Yx,\quad \text{ for all } x\in \fL^\perp,
\end{equation*}
proving the second equality in \eqref{relations} and the following inclusion
\begin{equation*}
\fL^\perp\subset \Ker(X-Y).
\end{equation*}
It remains to check the opposite inclusion
\begin{equation}\label{baba}
\fL^\perp\supset \Ker(X-Y).
\end{equation}
Let $x\in \Ker(X-Y)$ admit the representation
\begin{equation*}
x=\v+w,
\end{equation*}
where $\v\in \fL$ and $w\in \fL^\perp$. Then
\begin{equation}\label{dada}
 0=(X-Y)x=2X\v,
\end{equation}
using \eqref{relations}. Since $\fL\subset \Ker(I_{\fH_0}-X^\ast X)$, by \eqref{dada}
\begin{equation*}
0=\v-X^\ast X\v=\v,
\end{equation*}
that is, $x=w\in \fL^\perp$, proving \eqref{baba}. Thus,
\eqref{iadro} holds.

Finally, we prove \eqref{i4}. First, we notice that
\eqref{i3} and \eqref{Kontraktion:1} imply
that if $x\in\fL$ then $x\in\Ker(X+Y)$ and for any $y\in\Ker(X)$
\begin{equation*}
(x,y)=(X^\ast Xx,y)=(Xx,Xy)=0
\end{equation*}
Similarly, $(x,y)=0$ for any $y\in\Ker(Y)$. Hence, $x$ is orthogonal
to $\Ker(X)\cap\Ker(Y)$ and
\begin{equation}
\label{direct}
\fL\subset\Ker(X+Y)\ominus\big(\Ker(X)\cap\Ker(Y)\big).
\end{equation}

Suppose that the inverse inclusion does not hold. Then by
\eqref{direct} there is a nonzero
\begin{equation}
\label{yind}
y\in\Ker(X+Y)\ominus\big(\Ker(X)\cap\Ker(Y)\big)
\end{equation}
orthogonal to $\fL$, i.e. $y\in\fL^\perp$. By the second
equality in \eqref{relations} we will have $Xy=Yy$ which
contradicts \eqref{yind}. Hence equality \eqref{i4} holds true.

The proof is complete.
\end{proof}
Given any contractive solution $X$ associated with a \emph{spectral} subspace, Lemma
\ref{Neu:Theorem} allows one to provide a complete characterization of the set of all
contractive solutions to the Riccati equation in the sense that all contractive
solutions $Y$ to the Riccati equation
\eqref{Riccati} are in one-to-one correspondence with the closed
subspaces
\begin{equation}\label{nado}
\fL\subset\Ker(I_{\fH_0}-X^\ast X)\cap\Ker(XVX-V^\ast)
\end{equation}
reducing both the operators $A_0$ and $VX$. An explicit description of this correspondence is a
content of Theorem \ref{last:but:not:least} below. In particular, this theorem provides an
efficient criterion for a contractive solution $X$ associated with a spectral subspace of $\bB$ to
be unique.
\begin{theorem}\label{last:but:not:least}
Assume Hypothesis \ref{diag} and suppose that
$X\in\cB(\fH_0,\fH_1)$, $\|X\|\leq 1$ is a contractive solution
to the Riccati equation \eqref{Riccati} such that $\cG(\fH_0,X)$
is a spectral subspace of the operator $\bB$. Denote by $\cS$
the set of all contractive solutions to the Riccati equation
\eqref{Riccati} and by ${\cM}$ the lattice of all closed
subspaces of the Hilbert space $\fH_0$. Then the mapping
$\cT_X:\cS\longrightarrow {\cM}$ introduced by
\begin{equation*}
\cT_X(Y)=\Ker(I_{\fH_0}+Y^\ast X), \quad Y\in\cS
\end{equation*}
is one-to-one and the image of $\cT_X$ coincides
with the set $\cR $ of all closed subspaces
$\fL\subset \fH_0$ satisfying \eqref{nado} and
reducing both the operators $A_0$
and $VX$.

In particular, if $X\in\cB(\fH_0,\fH_1)$, is a contractive
solution to the Riccati equation \eqref{Riccati} associated with
a spectral subspace of the operator $\bB$, then $X$ is a unique
contractive solution to \eqref{Riccati} iff
\begin{equation}
\label{KKzero} \Ker(I_{\fH_0}-X^\ast X)\cap\Ker(XVX - V^\ast)=\{0\}.
\end{equation}
\end{theorem}
\begin{proof}
Let $Y\in\cS$ be arbitrary. Since the graph of $X$ is a spectral
subspace of $\bB$, the orthogonal projections onto the graphs of $X$ and $Y$ commute.
Then by Lemma \ref{Neu:Theorem}
\begin{equation}\label{nunu}
Y|_{\cT_X(Y)} = - X|_{\cT_X(Y)}\quad \text{and} \quad
Y|_{\cT_X(Y)^\perp}= X|_{\cT_X(Y)^\perp},
\end{equation}
which proves, in particular, that the mapping $\cT_X$ is one-to-one.

It remains to prove that
\begin{equation*}
\Ran\cT_X =\cR.
\end{equation*}
We start with the proof of the inclusion
\begin{equation}\label{ranll}
\Ran \cT_X\subset\cR.
\end{equation}
First, we prove that the subspace
\begin{equation}
\label{lll}
\fL=\cT_X(Y),\quad Y\in\cS
\end{equation}
reduces $A_0+VX$. That is, we need to establish that $\fL$ and $\fL^\perp$ are
$(A_0+VX)$-invariant subspaces.

The fact that $\fL^\perp$ is $(A_0+V X)$-invariant can be proven as follows. Taking
into account that both $X$ and $Y$ satisfy the Riccati equation \eqref{Riccati} and
by Lemma \ref{Neu:Theorem} $(X-Y)x=0$ for $x\in\fL^\perp$, a simple computation shows
that
\begin{equation*}
(X-Y)(A_0+VX)x=0
\quad\text{for any}\quad
x\in\fL^\perp.
\end{equation*}
Applying Lemma \ref{Neu:Theorem} again yields $(A_0+VX)x\in
\fL^\perp$ for any $x\in\fL^\perp$ which proves that $\fL^\perp$
is $(A_0+V X)$-invariant.

Next we establish that $\fL$ is $(A_0+V X)$-invariant. Since $\fL^\perp$ is $(A_0+V X)$-invariant,
the subspace $\fL$ is invariant for the operator $A_0+X^\ast V^\ast$. Note that the operator $(A_0
+ X^\ast V^\ast)(I_{\fH_0}+ X^* X)$ is self-adjoint. This fact is proven in
\cite{Langer:Tretter:2000b}, \cite{Motovilov:95} but alternatively can easily be seen from the
identity
\begin{equation*}
(x+Xx,\bB(x+Xx))=(x,(A_0 + X^\ast V^\ast)(I_{\fH_0}+X^\ast X)x) \quad\text{for any}\quad
x\in\fH_0.
\end{equation*}
Taking into account that by Lemma \ref{Neu:Theorem} $X^\ast X|_{\fL}=I_{\fL} $, one concludes that
\begin{equation*}
(I_{\fH_0} + X^\ast X)(A_0 + V X)x = (A_0 + X^\ast V^\ast)(I_{\fH_0}+X^\ast X)x =2(A_0 + X^\ast
V^\ast),\quad x\in \fL,
\end{equation*}
which implies $(A_0 + V X)x\in\fL$, proving that $\fL$ is also
$(A_0+VX)$-invariant. Thus we have proven that $\fL$ reduces the
operator $A_0+VX$.

The same arguments hold for the operator $A_0+V Y$. In particular, the subspace $\fL$
reduces the operator $A_0+VY$.

Now we are ready to prove inclusion \eqref{ranll}. Combining the facts that $\fL$
reduces $A_0+VX$ as well as $A_0+VY$ and that $X|_\fL=-Y|_\fL$ implies that $\fL$
reduces the operators $A_0$, $VX$, and $VY$. In particular,
\begin{align*}
0&=(A_1Y-YA_0-YVY+V^\ast)x\\
 &=(-A_1Xx+XA_0x-XVXx+V^\ast)x\\
 &=-2(XVX-V^\ast)x, \qquad x\in \fL,
\end{align*}
proving that
\begin{equation}
\label{pod} \fL\subset\Ker(XVX-V^\ast),
\end{equation}
and hence \eqref{nado} holds, since $\fL\subset \Ker(I_{\fH_0}-X^\ast X)$ by Lemma
\ref{Neu:Theorem}. Thus, $\fL=\cT_X(Y)\subset\cR$ which proves the inclusion
\eqref{ranll}.

In order to complete the proof of the theorem it remains to prove the
opposite inclusion
\begin{equation}
\label{ranl}
\cR\subset\Ran\cT_X.
\end{equation}
Let $\fL\subset\cR$ be arbitrary. Introduce the contraction $Y$ by
setting
\begin{equation}\label{lper}
Y|_{\fL} = -X|_{\fL},\quad \text{and} \quad Y|_{\fL^\perp} =
X|_{\fL^\perp}.
\end{equation}
We need to show that $Y\in\cS$ and that $\cT_X(Y)=\fL$.

For $x\in\fL^\perp$ one obtains
\begin{equation*}
(A_1Y-YA_0-YVY+V^\ast)x=(A_1X-XA_0-XVX+V^\ast)x = 0
\end{equation*}
using the invariance of $\fL^\perp$ with respect to the operators $A_0$ and $VX$, the
fact that $X$ solves the Riccati equation \eqref{Riccati}, and the second equality in
\eqref{lper}.

Using the invariance of $\fL$ with respect to the operators $A_0$ and $VX$, and the
first equality in \eqref{lper}, for $x\in\fL$ one obtains
\begin{equation*}
(A_1 Y - Y A_0 - Y V Y + V^\ast)x = \ (-A_1 X + XA_0 - XV X + V^\ast)x.
\end{equation*}
Since $\fL\subset\cR$, and hence $\fL\subset\Ker(XVX-V^\ast)$, for $x\in\fL$ one concludes that
$(XVX-V^\ast)x=0$. Therefore,
\begin{equation*}
(-A_1 X + XA_0 - XV X + V^\ast)x=(-A_1 X + XA_0 + XV X - V^\ast)x,
\end{equation*}
which is zero, since $X$ solves the Riccati equation \eqref{Riccati}.
Hence,
\begin{equation*}
(A_1 Y - Y A_0 - Y V Y + V^\ast)x = 0, \quad x\in \fL.
\end{equation*}
Thus, we constructed a contractive solution $Y$ to the Riccati
equation \eqref{Riccati}, which yields $Y\in\cS$. Applying
Lemma \ref{Neu:Theorem} implies that
\begin{equation}\label{lperr}
Y|_{\cT_X(Y)} = -X|_{\cT_X(Y)},\qquad Y|_{\cT_X(Y)^\perp} = X|_{\cT_X(Y)^\perp}.
\end{equation}
and
\begin{equation}\label{raz}
\cT_X(Y)\subset \Ker(I_{\fH_0}-X^\ast X).
\end{equation}
Since $\fL\subset\cR$ one also concludes that
\begin{equation}\label{dva}
\fL\subset \Ker(I_{\fH_0}-X^\ast X).
\end{equation}
Combining \eqref{lper}, \eqref{lperr}, \eqref{raz}, and
\eqref{dva} proves that
\begin{equation*}
\cT_X(Y)=\fL.
\end{equation*}
Thus, inclusion \eqref{ranl} is proven.

The proof is complete.
\end{proof}
\begin{remark}\label{ganz:neu0}
Notice that the subspace
\begin{equation*}
\fL'=\{x\oplus Xx|\ x\in\fL\},
\end{equation*}
where $\fL$ stands for any subspace referred to in Theorem
\ref{last:but:not:least}, is simultaneously $\bA$- and
$\bB$-invariant and, therefore, it can be split from the further
considerations if necessary.
\end{remark}
\begin{proof}
On the one hand,
\begin{align*}
\bA(x\oplus Xx) &= A_0x\oplus A_1 Xx\\
 &=A_0 x\oplus XA_0x, \quad x\in \fL
\end{align*}
since $x\in \fL\subset \Ker(XVX-V^\ast)$ taking into account that $X$ solves \eqref{Riccati},
proving that $\fL'$ is also $\bA$-invariant.

Since $X$ solves the Riccati equation \eqref{Riccati}, for any
$x\in\fH_0$, in particular, for $x\in \fL$ one has
\begin{align*}
\bB(x\oplus Xx) &= (A_0 + VX)x\oplus (V^\ast + A_1 X)x\\
&=(A_0 + VX)x\oplus X(A_0 + VX)x,
\end{align*}
which proves that $\fL'$ is also $\bB$-invariant, since $\fL$ is
$(A_0 + VX)$-invariant by hypothesis.
\end{proof}

As an immediate corollary of Theorem \ref{last:but:not:least} we
get the following uniqueness results.
\begin{corollary}\label{nagelneu} Assume Hypothesis \ref{diag}.
Let $X\in\cB(\fH_0,\fH_1)$, $\|X\|\leq 1$ be a contractive
solution to the Riccati equation \eqref{Riccati} such that
$\cG(\fH_0,X)$ is a spectral subspace of the operator $\bB$.

(i) If $X$ is a strictly contractive operator, i.e., $\|Xx\| <
\|x\|$ for any $x\in\fH_0$, $x\neq 0$, then $X$ is a unique
contractive solution to \eqref{Riccati}.

(ii) If
\begin{equation}
\label{dissi} \Ker(I_{\fH_0}-X^\ast X)\cap\Ker\bigl(\Im(VX)\bigr)=\{0\},
\end{equation}
then $X$ is a unique contractive solution to \eqref{Riccati}. In particular, if $VX$ is a
dissipative operator with positive imaginary part, them $X$ is a unique contractive solution to
\eqref{Riccati}.
\end{corollary}

\begin{proof}
(i) If $X$ is a strictly contractive operator, then $\Ker(I_{\fH_0}-X^\ast X)=\{0\}$. Hence
\eqref{KKzero} holds and by Theorem \ref{last:but:not:least} the operator $X$ is the unique
contractive solution to \eqref{Riccati}.

(ii) Suppose that  $Y$ is a contractive solution of the Riccati equation \eqref{Riccati}.
Introducing the subspace $\fL=\Ker(I_{\fH_0}+Y^\ast X)$, by Theorem \ref{last:but:not:least} one
concludes that
\begin{equation*}
\fL\subset\Ker(I_{\fH_0}-X^\ast X)\cap\Ker(XVX-V^\ast).
\end{equation*}
In particular,
\begin{equation*}
(XVX-V^\ast)x=0,\quad x\in \fL
\end{equation*}
and hence
\begin{equation}\label{XXXX}
X^\ast (XVX-V^\ast) x=X^\ast XVXx-X^\ast V^\ast x=0,\quad x\in \fL.
\end{equation}
By Theorem \ref{last:but:not:least} $\fL$ reduces the operator $VX$. In particular,
\begin{equation}
\label{XXVX} X^\ast XVXx=VXx, \quad x\in\fL.
\end{equation}
Combining \eqref{XXXX} and \eqref{XXVX} yields
\begin{equation*}
(VX-X^\ast V^\ast)x=0, \quad x\in \fL,
\end{equation*}
that is, $x=0$ for any $x\in \fL$, since $\Ker(VX-X^\ast V^\ast)=\{0\}$ by hypothesis. Hence
$\fL=\Ker(I_{\fH_0}+Y^\ast X)=\{0\}$ which proves that $Y=X$ using  Lemma \eqref{Neu:Theorem},
completing the proof.
\end{proof}
\begin{remark}
Statement (i) of Corollary \ref{nagelneu} concerning the spectral subspaces $\Ran$
$\EE_{\bB}(\Delta)$ associated with \emph{closed} Borel sets $\Delta$ of the real
axis appeared first in \cite{Adamyan:Langer:Tretter:2000a} with a somewhat different
strategy of the proof based on a description of maximal $J$-non-negative subspaces in
a Krein space.
\end{remark}
\begin{remark}
Some different uniqueness results for Riccati equations in finite-dimen\-sional
Hilbert spaces were obtained in \cite{Daughtry}. Note that the property for a
solution to the Riccati equation to be isolated is related to its stability
\cite{Campbell:Daugthry}. Stability of invariant subspaces is studied in
\cite{Adams:83}.
\end{remark}
To illustrate the statement of Theorem \ref{last:but:not:least}
suppose that $\fH=\fH_0\oplus\fH_1$ where $\fH_0$ and $\fH_1$
are copies of the same Hilbert space $\fK$, i.e.,
$\fH_0=\fH_1=\fK$. Assume that $A_0=A_1=0$ and $V=I_\fK$ is the
identity operator in $\fK$. Then the Riccati equation
\eqref{Riccati} (after the appropriate identification of the
copies $\fH_0$ and $\fH_1$) reads as $X^2 = I_\fK$ and it
obviously has the solution $X=I_{\fK}$ associated with the
eigenspace of $\bB$ corresponding to the eigenvalue one.
Therefore, $X$ is an isolated point in the set of all solutions.
Obviously,
\begin{equation*}
\Ker(I_\fK-X^\ast X)\cap\Ker(XVX - V^\ast)=\fK,
\end{equation*}
$VX=I_{\fK}$, and $\fK$ reduces $A_0=0$. Therefore, by Theorem \ref{uni} all solutions to the
Riccati equation can be uniquely parameterized by closed subspaces $\fL\subset\fK$. The case of
$\fL=\{0\}$ corresponds to the identity solution $X=I_\fK$, the case $\fL=\fK$ corresponds to the
solution $\widehat X=-X=-I_{\fK}$ which is also isolated being associated with the eigenspace of
$\bB$ corresponding to the eigenvalue negative one. If $\dim \fK>1$, all the other solutions to
the Riccati equation $X^2=I_\fK$ can be uniquely parameterized by the nontrivial subspaces
$\fL\subset \fK$ of nonzero codimension and thus correspond to invariant subspaces of the operator
$\bB$ which are not spectral ones. All those solutions $Y$ are unitary self-adjoint operators in
$\fK$ different from $I_\fK$ and $-I_\fK$, with $\fL=\Ker(Y+I_\fK)$, and hence
\begin{equation*}
\spec( \bB|_{\cG(\fH_0,Y)})=\{1, -1\}.
\end{equation*}

It is worth to note that given a solution $X$ to the Riccati equation \eqref{Riccati} (in contrast
to the hypothesis of Theorem \ref{last:but:not:least} not necessarily contractive and not
necessarily associated with a spectral subspace of the operator matrix $\bB$), the set $\cR$ of
all closed subspaces of $\fH_0$ reducing both the operators $A_0$ and $VX$ and satisfying
\eqref{nado} admits a dual description in terms of the corresponding subspaces of $\fH_1$.

In order to formulate the precise statement we introduce the set $\cR_\ast$ of all
closed subspaces $\fL_\ast\subset{\fH_1}$ reducing both $A_1$ and $V^\ast X^\ast$ and
satisfying
\begin{equation}\label{nado2}
\fL_*\subset\Ker(I_{\fH_1}-XX^\ast) \cap \Ker(X^\ast V^\ast X^\ast - V).
\end{equation}

\begin{proposition}
Let $X$ be a solution to the Riccati equation \eqref{Riccati}. Then, under the above notations,
the sets $\cR$ and $\cR_\ast$ are in one-to-one correspondence under the mapping
\begin{equation*}
\fL \mapsto X \fL, \quad \fL \in \cR.
\end{equation*}
In particular, the inverse mapping is given by $\fL_\ast \mapsto X^\ast \fL_\ast$,
$\fL_\ast\in\cR_\ast$.
\end{proposition}

\begin{proof}
Let $\fL\in\cR$ be arbitrary. Set $\fL_\ast = X \fL$. For any $x^\ast\in\fL_\ast$
there is a unique $x\in\fL$ such that $x_\ast=Xx$. To prove this, suppose to the
contrary that there is another element $y\neq x$ in $\fL$ such that $Xy=x_\ast$. Then
$X(x-y)=0$ and thus $X^\ast X(x-y)=0$. But $x,y\in\Ker(I_{\fH_0}-X^\ast X)$ and,
therefore, $X^\ast X(x-y)=x-y\neq0$. A contradiction.

If $x_\ast = Xx$, $x\in\fL$, then
\begin{equation*}
(A_1-V^\ast X^\ast)x_* = (A_1-V^\ast X^\ast)Xx = (A_1X-V^\ast)x = X(A_0-VX)x \in \fL_\ast
\end{equation*}
by successive use of \eqref{nado}, the hypothesis that $X$ solves the Riccati equation
\eqref{Riccati} and that $\fL$ is obviously $(A_0-VX)$-invariant. Thus $\fL_\ast$ is $(A_1-V^\ast
X^\ast)$-invariant. Moreover, $A_1 x_\ast=A_1X x=X A_0 x\in \fL_\ast$, since $x\in \fL \subset
\Ker(XVX-V^\ast)$ by \eqref{nado} and $X$ solves \eqref{Riccati}. Thus $\fL_\ast$ in addition is
$A_1$-invariant, proving that $\fL_\ast$ is invariant for both $A_1$ and $V^\ast X^\ast$. Further,
we note that $x_\ast=Xx\in \Ker(I_{\fH_1}-XX^\ast)$ and hence a simple computation
\begin{gather*}
(V-X^\ast V^\ast X^\ast)x_\ast = (V-X^\ast V^\ast X^\ast)Xx\\
=(X^\ast XV-X^\ast V^\ast X^\ast)Xx = X^\ast(XVX-V^\ast)x=0
\end{gather*}
proves inclusion \eqref{nado2}.

Now we claim that $\fL_\ast^\perp=\fH_1\ominus\fL_\ast$ is invariant for $V^\ast X^\ast$. To show
this choose arbitrary $y\in\fL_\ast^\perp$ and $x\in\fL_\ast$. Then
\begin{equation*}
(V^\ast X^\ast y, x) = (y, XVx) = (y, XX^\ast V^\ast X^\ast x)
\end{equation*}
since $Vx=X^\ast V^\ast X^\ast x$. The subspace $\fL_\ast$ is invariant for $V^\ast X^\ast$ and by
\eqref{nado2} we get
\begin{equation*}
X X^\ast V^\ast X^\ast x \in \fL_\ast.
\end{equation*}
Thus, $(V^\ast X^\ast y, x) = 0$. Since $x$ and $y$ are arbitrary, this implies that $V^\ast
X^\ast y\in \fL_\ast^\perp$.

Therefore, we proved that $\fL_\ast$ reduces both $A_1$ and $V^\ast X^\ast$. Thus,
$X\fL\in\cR^\ast$ and the mapping $\fL\mapsto X\fL$ maps $\cR$ onto $\cR_\ast$. By symmetry we
also conclude that this mapping is one-to-one.
\end{proof}


\end{document}